\documentclass[12pt,naturalnames]{article}
\date{October 21, 2002}
\title%
{First Passage Percolation Has\\Sublinear Distance Variance}
\author{ Itai Benjamini \and Gil Kalai \and Oded Schramm }

\newif\ifdraft
\drafttrue
\newif\ifhyper
\hypertrue  
\usepackage{amsmath}
\usepackage{amsthm}
\usepackage{amsfonts}
\newif\ifhyper\IfFileExists{hyperref.sty}{\hypertrue}{\hyperfalse}
\ifhyper\usepackage{hyperref}\fi

\theoremstyle{plain}
\newtheorem{theo}{Theorem}
\newtheorem{thm}[theo]{Theorem}
\newtheorem{theorem}[theo]{Theorem}
\newtheorem{lemma}[theo]{Lemma}

\theoremstyle{remark}

\def\Bbb#1{\mathbb{#1}}
\newcommand{\R}{{\Bbb R}}
\newcommand{\Z}{{\Bbb Z}}

\newcommand{\N}{{\Bbb N}}
\DeclareMathOperator{\length}{length}

\def\P{\PP}
\def\PP{{\mathbf P}}

\def\E{{\bf E}}

\def\Bb#1#2{{\def\md{\bigm| }#1\bigl[\,#2\,\bigr]}}
\def\BB#1#2{{\def\md{\Bigm| }#1\Bigl[\,#2\,\Bigr]}}
\def\Bs#1#2{{\def\md{\mid}#1[\,#2\,]}}
\def\Pb{\Bb\PP}
\def\PB{\BB\P}
\def\Ps{\Bs\PP}
\def\Eb{\Bb\E}

\def\Es{\Bs\E}

\def\md{\mid}
\def\eref#1{(\ref{#1})}

\def\Zd{\Z^d}
\def\EE{\mathcal E}
\def\EE{{E}}

\def\le{\leq}
\def\ge{\geq}
\def\hat{\widehat}
\def\var{{\mathrm {var}}}

\def\dist{{\mathrm{dist}}}

\let\CHI=\chi
\def\chi{\raise2.2pt\hbox{$\CHI$}}

\def\proofstyle#1{\medskip\noindent{\sc #1: }}
\def\proof{\proofstyle{Proof}}

\def\QED{\qed\medskip}

\long\def\comment#1{}

\begin{document}

\maketitle

\begin {abstract}
Let $0<a<b<\infty$, and for each edge $e$ of $\Z^d$ let
$\omega_e=a$ or $\omega_e=b$, each with probability $1/2$,
independently.  This induces a random metric $\dist_\omega$
on the vertices of $\Z^d$, called first passage percolation. 
We prove that for $d>1$ the distance $\dist_\omega(0,v)$ from the origin
to a vertex $v$, $|v|>2$, has variance bounded by $C\,|v|/\log|v|$,
where $C=C(a,b,d)$ is a constant which may only depend on $a$, $b$
and $d$.  Some related variants are also discussed.
\end {abstract}


\section{Introduction}

Consider the following model of first passage percolation.
Fix some $d=2,3,\dots$, and let $\EE=\EE(\Z^d)$ denote the set
of edges in $\Z^d$.  Also fix numbers $0<a<b<\infty$.
Let $\Omega:=\{a,b\}^\EE$ carry the product measure,
where $\Ps{\omega_e=a}=\Ps{\omega_e=b}=1/2$ for each $e\in\EE$.
Given $\omega=(\omega_e:e\in E)\in\Omega$ and vertices $v,u\in\Zd$,
let $\dist_\omega(v,u)$ denote the least distance from
$v$ to $u$ in the metric induced by $\omega$; that is,
the infimum of $\length_\omega(\alpha):=\sum_{e\in\alpha}\omega_e$,
where $\alpha$ ranges over all finite paths in $\Zd$ from $u$ to $v$. 
Let $|v|:=\|v\|_1$ for vertices $v\in\Zd$.

\begin{thm}\label{t.fpp}
There is a constant $C=C(d,a,b)$ such that for every $v\in\Zd$, $|v|\ge 2$,
$$
\var\bigl(\dist_\omega(0,v)\bigr)
\le
C\,\frac{|v|}{ \log|v|}\,.
$$
\end{thm}

In  \cite {Ke2} Kesten used
martingale inequalities to prove 
$\var\bigl(\dist_\omega(0,v)\bigr) \le C\,|v|$
and proved some tail estimates. Talagrand \cite {T} used his 
``convexified'' discrete isoperimetric inequality to prove 
that  for all $t>0$,
\begin {equation}
\label {talatala}
\PB { |\dist_\omega(0,v) -M| \ge t \sqrt {|v|} } \le C\, \exp (-t^2/C)\,,
\end {equation}
where
$M$ is the median value of $\dist_\omega(0,v)$.
(Both Kesten's and Talagrand's results apply to more general distributions
of the edge lengths $\omega_e$.)
``Novice readers might expect to hear next of a central limit theorem
being proved,'' writes Durrett \cite {Du}, describing Kesten's results,
``however physicists tell us \dots\ that in two dimensions the standard deviation
\dots\ is of order $|v|^{1/3}$.'' 
Recent remarkable work
\cite {BDJ,J,DZ} not only supports this prediction, but suggests
what the limiting distribution and large deviation behavior is.
The case of a certain
variant of {\bf oriented} first passage percolation
is settled by Johansson \cite {J}.
For lower bounds on the variance in two-dimensional
first passage percolation see
Newman-Piza~\cite{NP} and Pemantle-Peres~\cite{PP}.

\medskip

As in Kesten's and Talagrand's earlier results, the most
essential feature about first passage percolation which
we use is that the number of edges $e\in \EE$
such that modifying $\omega_e$ increases $\dist_\omega(0,v)$
is bounded by $C\,|v|$. 

Another essential ingredient in the current paper is the following
extension by Talagrand~\cite{T1} of an inequality by
Kahn, Kalai and Linial \cite {KKL}.
Let $J$ be a finite index set.
For $j\in J$, and $\omega\in\{a,b\}^J$ let $\sigma_j\,\omega$ be
the element of $\{a,b\}^J$ which is different from $\omega$
only in the $j$-th coordinate.
For $f:\{a,b\}^J\to\R$ set $\sigma_j f:=f\circ\sigma_j$ and
$$
\rho_jf(\omega):=
\frac{f(\omega)-\sigma_jf(\omega)}{2}\,.
$$
Talagrand's~\cite[Thm.~1.5]{T1} inequality is 
\begin{equation}\label{e.tala}
\var (f) \le C\, \sum_{j\in J} \frac{\|\rho_jf\|^2_2}
{1+\log\bigl (\|\rho_jf\|_2/\|\rho_jf\|_1\bigr)}\,,
\end{equation}
where $C$ is a universal constant.  (In Section~\ref{s.tala}
we supply a direct proof of~\eref{e.tala} from
the Bonami-Beckner inequality.  
A very reasonable upper bound for $C$ can be obtained from this
proof.  We also explain there how~\eref{e.tala} easily
implies~\eref{talatala}.)

The basic idea in the proof of Theorem~\ref{t.fpp} is to apply this
inequality to $f(\omega):=\dist_\omega(0,v)$.  
For this, we wanted to show that, roughly,
$
\Pb{\rho_e f(\omega)\ne 0}
$
is small, except for a small number of edges $e$.
However, since we were not able to prove this,
we had to resort to an averaging trick.

\medskip

Theorem~\ref{t.fpp} should hold for other models
of first passage percolation, where the edge lengths
have more general distributions. 
The relevant result of~\cite{KKL}, as well as~\eref{e.tala} that we use
rely on the Bonami-Beckner~\cite{Bo,Be} inequality,
which holds for $\{a,b\}^n$, but fails on some more general product spaces.
The paper \cite {BKKKL} does extend some of these results to
general product spaces; see also~\cite{Fr-kkl}.
Talagrand's~\cite[Thm.~1.5]{T1} applies to product measures on
$\{a,b\}^n$, which are not necessarily uniform.
Beckner~\cite{Be} uses his inequality to derive a 
similar result for the Gaussian measure on
$\R^n$, and an analog for Talagrand's inequality~\eref{e.tala}  
for the Gaussian measure
(pointed out to us by Assaf Naor) can be found in~\cite {BH}.
In the present paper, we preferred simplicity to generality,
but it will be interesting to extend the theorem
to more general distributions. We refer the reader also to the book
by Ledoux~\cite {Ledoux} for a 
general view of relevant techniques and knowledge.

It seems even more important to put effort into the fundamental
task of lowering the upper bound from $C\,|v|/\log|v|$
to $C\,|v|^{1-\epsilon}$, $\epsilon>0$.

\medskip

To first illustrate the basic technique in a simpler setting,
we will also prove a theorem about
the variance of the first passage percolation circumference
of a discrete torus, or more generally in a cartesian product of
a finite vertex-transitive graph with a circle.

Let $H$ be a finite vertex-transitive graph, and let $n\in\N_+$.
Let $G$ be the product of $H$ and the cycle of length $n$, $\Z/(n\,\Z)$.
We say that a closed path $\gamma$ in $G$ is a circumference, if its projection onto
$\Z/(n\,\Z)$ has degree $1$; namely, $\gamma$ can be oriented so that its projection
has one more edge going from $0$ to $1$ than from $1$ to $0$.
If $\omega:E(G)\to\{a,b\}$, define its circumference length $c_G(\omega)$
as the minimal $\omega$-length of any circumference path.

\begin{theorem}\label{t.circum}
Let $n\in\N_+$ and let $G$ be a cartesian product of a vertex-transitive graph $H$
with the cycle $\Z/(n\,\Z)$ of length $n$. 
Let $\omega\in\{a,b\}^{E(G)}$ be random-uniform.
Then
$$
\var\bigl(c_G(\omega)\bigr)
\le C\,\frac ba(b-a)^2\,\frac {n}{1+\log(a\,|V(H)|/b)}
\,,
$$
where $C$ is a universal constant.
\end{theorem}

For example, when $G$ is the square torus $G=\bigl(\Z/(n\,\Z)\bigr)^2$, with $n>1$,
we get the estimate $\var(c_G)\le C\, (b/a)\,(b-a)^2\,n/\log n$.

\section{Proof of Theorem~\ref{t.circum}}

Let $\beta$ be a circumference path in $G$
such that $\sum_{e\in\beta}\omega_e =c_G(\omega)$.
We use some arbitrary but fixed method for choosing
between all possible choices for $\beta$.
Clearly, $c_G(\omega)\le b\,n$.
Therefore,
\begin{equation}\label{e.beta}
|\beta|\le
b\,n/a\,,
\end{equation}
where $|\beta|$ denotes the number of edges in $\beta$.

Let $e\in E(G)$.
Note that if $\rho_e\,c_G(\omega)<0$, then we must have $e\in\beta$.
By~\eref{e.beta}, this gives
\begin{equation}\label{e.esum}
\sum_{e\in E(G)}\Ps{\rho_e\,c_G<0} \le
\Eb{|\beta|}\le b\,n/a\,.
\end{equation}
Let $\Gamma$ be the automorphism group of $G$.
Fix $e\in E(G)$.
By symmetry, $\Ps{\rho_e c_g<0}=\Ps{\rho_{e'}c_g<0}$ for every
$e'\in\Gamma e$.
Consequently, 
\begin{equation*}
|\Gamma e|\,\Ps{\rho_e c_G<0}\le
\sum_{e'\in E(G)}\Ps{\rho_{e'} c_G<0}\le b\,n/a
\,.
\end{equation*}
Now note that since $H$ is vertex-transitive, also $G$ is vertex-transitive,
and consequently $|\Gamma e|\ge |V(G)|/2=n\,|V(H)|/2$.
Therefore, 
\begin{equation}
\Ps{\rho_e c_G<0}\le b\,\bigl|V(H)\bigr|^{-1}/a
\,.
\label{e.den1}
\end{equation}
Now clearly, $\Ps{\rho_e c_G\ne 0}=2\,\Ps{\rho_e c_G<0}$
and $\|\rho_e c_G\|_\infty\le (b-a)/2$.  Therefore,
\begin{equation}\label{e.es}
\|\rho_e c_G\|_2^2 \le (1/2)\,(b-a)^2\,\Ps{\rho_e c_G<0}\,.
\end{equation}
By Cauchy-Schwarz,
\begin{equation}\label{e.den2}
\|\rho_e c_G\|_1 \le \sqrt{\Ps{\rho_e c_G\ne 0}}\,\|\rho_e c_G\|_2\,.
\end{equation}
By~\eref{e.tala}, we have
\begin{equation*}
\var(c_G)
\le C\,\frac {\sum_{e\in E(G)} \|\rho_e c_G\|_2^2}
{1+\min_{e\in E(G)} \log \bigl(\|\rho_e\, c_G\|_2/\|\rho_e\, c_G\|_1\bigr)}\,.
\end{equation*}
To estimate the numerator, we use~\eref{e.es}
and~\eref{e.esum}, and for the
denominator~\eref{e.den1} and~\eref{e.den2}.  The theorem
easily follows.
\QED

\section{Proof of Theorem~\ref{t.fpp}}\label{s.fpp}

If $v,u\in\Z^d$, and $\alpha$ is a path from $v$ to $u$,
then $\alpha$ will be called an $\omega$-{\bf geodesic} if
it minimizes $\omega$-length; that is,
$\dist_\omega=\sum_{e\in\alpha}\omega_e$.
Given $\omega$, let $\gamma$ be an $\omega$-geodesic
from $0$ to $v$. 
Although there may be more than one such geodesic,
we require that $\gamma$ depends only on $\omega$.
(For example, we may use an arbitrary deterministic choice
among any possible collection of $\omega$-geodesics.)

The general strategy for the proof of Theorem~\ref{t.fpp} is as for
Theorem~\ref{t.diamfpp}.  However, the difficulty is that there
is not enough symmetry to get a good bound on $\Pb{\rho_e\, \dist_\omega(0,v)<0}$.
It would have been enough to show that
$ \Pb{e\in\gamma}<C\,|v|^{-1/C} $
holds with the exception of at most $C\,|v|/\log|v|$ edges, for some
constant $C>0$.  But we could not prove this.
Therefore, we will need an averaging argument, for which the following
lemma will be useful.

\begin{lemma}\label{l.slow}
There is a constant $c>0$ such that for every $m\in\N$ there
is a function
$$
g=g_m:\{a,b\}^{m^2}\to \{0,1,2,\dots,m\}
$$
satisfying 
$$
\|\sigma_j\, g-g\|_\infty\le 1
$$
for every $j=1,2,\dots,m^2$
and 
$$
\max_y\Ps{g(x)=y}\le c/m\,,
$$
where $x$ is random-uniform in $\{a,b\}^{m^2}$.
\end{lemma}

\proof
Assume $a=0$ and $b=1$, for simplicity of notation.
Let $k:\N\to\{0,1,2,\dots,m\}$ be the
function satisfying $k(0)=0$,
$k(j+1)=k(j)+1$ when $j\in \bigcup_{s=0}^\infty [2sm,2sm+m-1]$
and $k(j+1)=k(j)-1$ for all other $j\in\N$.
It is left to the reader to check that
$g(x)=k(\|x\|_1)$ has the required properties.
\QED

\medskip

Fix some $v\in\Zd$ with $|v|$ large,
and set $f=f(\omega):= \dist_\omega(0,v)$.
Clearly, 
$f(\omega)\le b\,|v|$, and therefore $|\gamma|\le (b/a)\,|v|$,
where $|\gamma|$ denotes the number of edges in
$\gamma$.  In particular, $f$ depends on only
finitely many of the coordinates in $\omega$.
Also, $|\gamma|\le (b/a)|v|$ implies
\begin{equation}
\label{e.nu}
\sum_{e\in\EE}\Ps{e\in\gamma}=
\Eb{|\gamma|}\le (b/a)\,|v|
\,.
\end{equation}

Fix $m:=\lfloor |v|^{1/4}\rfloor$, and
let $S:= {\{1,\dots,d\}\times\{1,\dots,m^2\}}$.
Let $c>0$ and $g=g_m$ be as in Lemma~\ref{l.slow}.
Given $x=\bigl(x_{i,j}:\{i,j\}\in S\bigr)\in\{0,1\}^{S}$
let 
$$
z=z(x):=
\sum_{i=1}^d\,g(x_{i,1},\dots,x_{i,m^2})\, \mathbf e_i\,,
$$
where $\{\mathbf e_1,\dots,\mathbf e_d\}$ is the standard basis for $\R^d$.
Define 
$$
\tilde f(x,\omega):= \dist_\omega(z,v+z)\,.
$$
We think of $\tilde f$ as a function on the
space $\tilde\Omega:=\{a,b\}^{S\cup E}$.
Since $|z|\le md$, it follows that $|\tilde f-f|\le 2mdb$.
In particular,
\begin{equation}\label{e.cvar}
\var(f)\le \var (\tilde f)+4mdb\sqrt{\var(\tilde f)}+4 m^2d^2b^2\,.
\end{equation}
It therefore suffices to find a good estimate for $\var(\tilde f)$.

Let $e\in E$ be some edge.
We want to estimate its influence:
\begin{equation*}
I_e(\tilde f) := \Pb{ \sigma_e\tilde f(x,\omega)\ne\tilde f(x,\omega)}
=
2\,\Pb{ \sigma_e\tilde f(x,\omega)>\tilde f(x,\omega)}\,.
\end{equation*}
(Here, $\P$ is the uniform measure on $\tilde\Omega$.)
Note that if the pair $(x,\omega)\in\tilde\Omega$ satisfies
$\sigma_e\tilde f(x,\omega)>\tilde f(x,\omega)$,
then $e$ must be on every $\omega$-geodesic from
$z$ to $v+z$.
Consequently, conditioning on $z$ and translating $\omega$ and $e$ by $-z$ gives
\begin{equation}\label{e.ietf}
I_e(\tilde f) = 2\,\Pb{ \sigma_e\tilde f(x,\omega)>\tilde f(x,\omega)}
\le 2\, \Pb{e-z \in\gamma}\,.
\end{equation}
Let $Q$ be the set of edges $e'\in E(\Z^d)$ such that
$\Ps{e-z=e'}>0$.  The $L^1$ diameter
of $Q$ is $O(m)$.
(We allow the constants in the $O(\cdot)$ notation to depend
on $d,a$ and $b$, but not on $v$.)
Hence, the diameter of $Q$ in
the $\dist_\omega$ metric is also $O(m)$,
and therefore $|\gamma\cap Q|\le O(m)$.
But the lemma gives
$$
\max_{z_0}\Pb{z=z_0\md \omega}\le (c/m)^d\,.
$$
By conditioning on $\gamma$ and summing over the edges in
$\gamma\cap Q$, we therefore get
$$
\Pb{e\in\gamma+z\md\gamma}\le O(1)\, m^{1-d}\,.
$$
Consequently,~\eref{e.ietf} and the choice of $m$ give
\begin{equation}\label{e.ie}
I_e(\tilde f) 
\le
O(1)\,|v|^{-1/4}
\,.
\end{equation}
Also,~\eref{e.nu} implies
$$
\sum_{e\in E}\Pb{e-z \in\gamma\md z} \le
 (b/a)\,|v|
\,.
$$
Combining this with~\eref{e.ietf} therefore gives
$$
\sum_{e\in E}I_e(\tilde f) \le 2\, (b/a)\,|v|\,.
$$
Applying~\eref{e.ie} yields
\begin{equation}\label{e.iesum}
\sum_{e\in E}\frac{ I_e(\tilde f)}{1+\bigl|\log I_e(\tilde f)\bigr|}
\le
O(1)\,|v|/\log|v|
\,.
\end{equation}
On the other hand, $I_s(\tilde f)\le 2(b-a)$ for $s\in S$.
As $|S|= O(1)\,|v|^{1/2}$ and
$\|\rho_q f\|_\infty=O(1)$ for $q\in S\cup E$,
we get from~\eref{e.iesum} and~\eref{e.tala}
$$
\var(\tilde f)
\le
O(1)\,
\sum_{q\in E\cup S}\frac{ I_q(\tilde f)}{1+\bigl|\log I_q(\tilde f)\bigr|}
\le
O(1)\,|v|/\log|v|
\,.
$$
Therefore, Theorem~\ref{t.fpp} now follows from \eref{e.cvar}.
\QED

As alluded to in the introduction, the proof would have been simpler
if we could show that there is
a $C>0$ such that $ \Pb{\rho_e f(\omega)\ne 0}<|v|^{-C}$ holds with
the possible exception of $|v|/\log|v|$ edges.  
It would be interesting to prove the closely related statement that the
probability that $\gamma$ passes within distance $1$ from $v/2$
tends to zero as $|v|\to\infty$.

\section{A proof of Talagrand's inequality~\eref{e.tala}}\label{s.tala}

To prove~\eref{e.tala}, it clearly suffices to take $a=0,b=1$.
For $f:\{0,1\}^J\to\R$,
consider the Fourier-Walsh expansion of $f$,
$$ f = \sum_{S \subset J} \hat{f}(S)\, u_S\,,$$
where $u_S(\omega)=(-1)^{S\cdot \omega}$ and
$S\cdot \omega$ is shorthand for $\sum_{s\in S} \omega_s$. 
For each $p\in\R$ define the operator
\begin {equation*}
T_p (f):=\sum_{S\subset J}p^{|S|} \hat f(S) u_S\,,
\end {equation*}
which is of central importance in harmonic analysis.
The Bonami-Beckner~\cite{Bo,Be} inequality asserts that
\begin{equation}
\label{e.bb}
\|T_pf\|_2\le \|f\|_{1+p^2}
\,.
\end{equation}

Set $f_j:=\rho_j\, f$.
Because $\rho_j\,u_S=u_S$ if $j\in S$ and
$\rho_j\,u_S=0$ if $j\notin S$,
we have
$$
\hat{f_j}(S)=\begin{cases} \hat f(S)& j \in S\,,\\
 0&j\notin S\,.
\end{cases}
$$
Since $\|g\|_2^2=\sum_{S\subset J}\hat g(S)^2$,
it follows that 
$$
\int_0^1 \|T_p\, f_j\|_2^2\,dp
=
\int_0^1\sum_{S\subset J}p^{2|S|} \hat {f_j}(S)^2\,dp
=
\sum_{S\subset J}1_{j\in S}\,\frac{\hat f(S)^2}{2|S|+1}\,.
$$
Consequently,
$$
\sum_{j\in J}
\int_0^1 \|T_p\, f_j\|_2^2\,dp
=
\sum_{S\subset J}|S|\,\frac{\hat f(S)^2}{2|S|+1}
\ge \frac13\sum_{\emptyset \ne S\subset J}\hat f(S)^2=\frac{\var(f)}3
\,.
$$
Therefore,~\eref{e.bb} gives
\begin{equation}\label{e.fin}
\var(f)
\le
3
\sum_{j\in J} \int_0^1 \|f_j\|_{1+p^2}^2\,dp
 \,.
\end{equation}
An instance of the H\"older inequality
$$
\Es{|f_j|^{1+p^2}}\le \Es{f_j^2}^{p^2}\,\Es{|f_j|}^{1-p^2}
$$
implies
\begin{align*}
\int_{0}^1 \|f_j\|_{1+p^2}^2\,dp
&
\le
\int_{0}^1
\Bigl(\Eb{f_j^2}^{p^2}\,\Eb{|f_j|}^{1-p^2}\Bigr)^{2/(1+p^2)}
\,dp
\\ &
=
\|f_j\|_2^2\int_{0}^1
\bigl(\|f_j\|_1/\|f_j\|_2\bigr)^{2(1-p^2)/(1+p^2)}
\,dp 
\\ &
\le 
2\,
\|f_j\|_2^2
\int_{1/2}^1
\bigl(\|f_j\|_1/\|f_j\|_2\bigr)^{2(1-p^2)/(1+p^2)}
\,dp 
\,.
\end{align*}
Let $s(p):= 2(1-p^2)/(1+p^2)$. 
Since $s'(p)\le s'(1)=-2$ when $p\in [1/2,1]$,
the above gives
\begin{align*}
\int_{0}^1 \|f_j\|_{1+p^2}^2\,dp
&
\le 
2\,\|f_j\|_2^2
\int_{s(1/2)}^{s(1)}
\bigl(\|f_j\|_1/\|f_j\|_2\bigr)^{s} 
\,\frac{ds }{s'(p)}
\\&
\le
\|f_j\|_2^2
\int_{0}^{6/5}
\bigl(\|f_j\|_1/\|f_j\|_2\bigr)^{s} 
\,ds 
\\&
=
\|f_j\|_2^2
\,\frac 
{1- \bigl(\|f_j\|_1/\|f_j\|_2\bigr)^{6/5} }
{\log \bigl(\|f_j\|_2/\|f_j\|_1\bigr) }
\,.
\end{align*}
Now~\eref{e.fin} implies
\begin{equation}\label{e.ourtala} 
\var(f)\le 3\sum_{j\in J}
\|f_j\|_2^2
\,\frac 
{1- \bigl(\|f_j\|_1/\|f_j\|_2\bigr)^{6/5} }
{\log \bigl(\|f_j\|_2/\|f_j\|_1\bigr) }
\end{equation}
from which~\eref{e.tala} follows.
\QED

\medskip\noindent
{\bf Remark:} It is worth noting that the tail estimate
(\ref {talatala}) can 
be derived from the variance inequality (\ref {e.tala}). 
Indeed, let $f=\dist_\omega (0,v)$ and $u\in(0,3/4)$.
Set $s=s(u):=\inf\bigl\{t\ge 0:\Ps{f(\omega)>t}<u\bigr\}$,
and define $g(\omega):= \max\bigl \{f(\omega),s\bigr\}$.
It follows from (\ref {e.tala}) that 
$$\var (g) \le C\,u\, |v| /
\bigl ( 1+\min_e \log(\|\rho_e\,g\|_2/\|\rho_e\,g\|_1)  \bigr )
\,,$$
and hence $\var (g) \le C \,u\, |v| /\log (1/u)$. 
Therefore, there is a constant $C'$ such that
 $$\PB {f(x)> s(u)+C'\sqrt {|v|  /\log (1/u) }}< u/2.$$ 
That is, $s(u/2)\le s(u)+C'\sqrt {|v|  /\log (1/u) }$.
Induction gives for $n=1,2,\dots$
$$
s(2^{-n})\le s(1/2) + O(1)\,\sqrt {n\,|v|}\,,
$$
which is the ``upper tail'' bound from (\ref {talatala}). 
The lower tail is treated similarly.
(This proof is fairly general. 
Using the more specific arguments from Section~\ref{s.fpp},
a slight improvement for 
the tail estimates in certain ranges may be obtained.)

\bigskip\noindent{\bf Acknowledgements.}
We are most grateful to Elchanan Mossel for very useful
discussions and to Jan Vondr\'ak for detecting a mistake in what
was Theorem 2 in a previous version of the paper.

\bigskip
\filbreak
{
\small
\begin{sc}
\parindent=0pt\baselineskip=12pt
\def\email#1{\par\qquad {\tt #1} \smallskip}
\def\emailwww#1#2{\par\qquad {\tt #1}\par\qquad {\tt #2}\medskip}

The Weizmann Institute of Science,
Rehovot 76100, Israel
\emailwww{itai@wisdom.weizmann.ac.il}
{http://www.wisdom.weizmann.ac.il/\string~itai/}

The Hebrew University,
Givat Ram, Jerusalem 91904,
Israel
\emailwww{kalai@math.huji.ac.il}
{http://www.ma.huji.ac.il/\string~kalai/}

Microsoft Research,
One Microsoft Way,
Redmond, WA 98052,
USA
\emailwww{schramm@microsoft.com}
{http://research.microsoft.com/\string~schramm/}
\end{sc}
}
\filbreak

\end{document}